\documentclass[12pt]{amsart}
\usepackage{amssymb}
\usepackage{color}
\usepackage{amsmath,epic,curves,amscd}
\usepackage[english]{babel}
\usepackage{graphicx}
\usepackage{comment}
\usepackage{appendix}
\usepackage{mathdots}
\usepackage[all]{xy}
\usepackage{mathtools}
\usepackage{lscape}
\usepackage{stmaryrd}
\usepackage{mathabx}
\theoremstyle{remark}

\allowdisplaybreaks \textwidth 15.5 true cm \textheight 24.2 true cm
\usepackage{color}

\begin{document}
\baselineskip 6.0 truemm
\parindent 1.5 true pc

\newcommand\lan{\langle}
\newcommand\ran{\rangle}
\newcommand\tr{\operatorname{Tr}}
\newcommand\ot{\otimes}
\newcommand\ttt{{\text{\rm t}}}
\newcommand\rank{\ {\text{\rm rank of}}\ }
\newcommand\choi{{\rm C}}
\newcommand\dual{\star}
\newcommand\flip{\star}
\newcommand\cp{{{\mathbb C}{\mathbb P}}}
\newcommand\ccp{{{\mathbb C}{\mathbb C}{\mathbb P}}}
\newcommand\pos{{\mathcal P}}
\newcommand\tcone{T}
\newcommand\mcone{K}
\newcommand\superpos{{{\mathbb S\mathbb P}}}
\newcommand\blockpos{{{\mathcal B\mathcal P}}}
\newcommand\jc{{\text{\rm JC}}}
\newcommand\dec{{\mathbb D}{\mathbb E}{\mathbb C}}
\newcommand\decmat{{\mathcal D}{\mathcal E}{\mathcal C}}
\newcommand\ppt{{\mathcal P}{\mathcal P}{\mathcal T}}
\newcommand\pptmap{{\mathbb P}{\mathbb P}{\mathbb T}}
\newcommand\xxxx{\bigskip\par ================================}
\newcommand\join{\vee}
\newcommand\meet{\wedge}
\newcommand\ad{\operatorname{Ad}}
\newcommand\ldual{\varolessthan}
\newcommand\rdual{\varogreaterthan}
\newcommand{\slmp}{{\mathcal M}^{\text{\rm L}}}
\newcommand{\srmp}{{\mathcal M}^{\text{\rm R}}}
\newcommand{\smp}{{\mathcal M}}
\newcommand{\id}{{\text{\rm id}}}
\newcommand\tsum{\textstyle\sum}
\newcommand\hada{\Theta}
\newcommand\ampl{\mathbb A^{\text{\rm L}}}
\newcommand\ampr{\mathbb A^{\text{\rm R}}}
\newcommand\amp{\mathbb A}
\newcommand\rk{{\text{\rm rank}}\,}
\newcommand\calI{{\mathcal I}}
\newcommand\bfi{{\bf i}}
\newcommand\bfj{{\bf j}}
\newcommand\bfk{{\bf k}}
\newcommand\bfl{{\bf l}}
\newcommand\bfzero{{\bf 0}}
\newcommand\bfone{{\bf 1}}
\newcommand\calc{{\mathcal C}}
\newcommand\calm{{\mathcal M}}
\newcommand\calb{{\mathcal B}}

\title{Choi matrices revisited}

\author{Seung-Hyeok Kye}
\address{Department of Mathematics and Institute of Mathematics, Seoul National University, Seoul 151-742, Korea}
\email{kye at snu.ac.kr}

\thanks{partially supported by NRF-2020R1A2C1A01004587, Korea}

\begin{abstract}
A linear map between matrix algebras corresponds to the Choi matrix
in the tensor product of two matrix algebras,
whose definition depends on the matrix units.
Paulsen and Shultz [J. Math. Phys. {\bf 54} (2013), 072201] considered the question
if one can replace matrix units by another basis of matrix algebras in the definition of Choi matrix to retain the correspondence
between complete positivity of maps and positivity of Choi matrices, and gave a sufficient condition
on basis under which this is true. In this note, we provide necessary and sufficient conditions, to see that
the Paulsen--Shultz condition is also necessary.
\end{abstract}
\maketitle

%%%%%%%%%%%%%%%%%%%%%%%%%%%%%%%%%%%%%%%%%%%%%%%%%%%%%%%%%%%%%%%%%%%%%%%%%%%%%%%%%%%%%%%%%%%%%%%%%%%%%%%%%%%%%%%%%%%%%%%%%%%%%%%%
%%%%%%%%%%%%%%%%%%%%%%%%%%%%%%%%%%%%%%%%%%%%%%%%%%%%%%%%%%%%%%%%%%%%%%%%%%%%%%%%%%%%%%%%%%%%%%%%%%%%%%%%%%%%%%%%%%%%%%%%%%%%%%%%
%%%%%%%%%%%%%%%%%%%%%%%%%%%%%%%%%%%%%%%%%%%%%%%%%%%%%%%%%%%%%%%%%%%%%%%%%%%%%%%%%%%%%%%%%%%%%%%%%%%%%%%%%%%%%%%%%%%%%%%%%%%%%%%%
%%%%%%%%%%%%%%%%%%%%%%%%%%%%%%%%%%%%%%%%%%%%%%%%%%%%%%%%%%%%%%%%%%%%%%%%%%%%%%%%%%%%%%%%%%%%%%%%%%%%%%%%%%%%%%%%%%%%%%%%%%%%%%%%
%%%%%%%%%%%%%%%%%%%%%%%%%%%%%%%%%%%%%%%%%%%%%%%%%%%%%%%%%%%%%%%%%%%%%%%%%%%%%%%%%%%%%%%%%%%%%%%%%%%%%%%%%%%%%%%%%%%%%%%%%%%%%%%%
%%%%%%%%%%%%%%%%%%%%%%%%%%%%%%%%%%%%%%%%%%%%%%%%%%%%%%%%%%%%%%%%%%%%%%%%%%%%%%%%%%%%%%%%%%%%%%%%%%%%%%%%%%%%%%%%%%%%%%%%%%%%%%%%
%%%%%%%%%%%%%%%%%%%%%%%%%%%%%%%%%%%%%%%%%%%%%%%%%%%%%%%%%%%%%%%%%%%%%%%%%%%%%%%%%%%%%%%%%%%%%%%%%%%%%%%%%%%%%%%%%%%%%%%%%%%%%%%%
%\section{Introduction}

For a given linear map $\phi$ from the algebra $M_n$ of all $n\times n$ complex matrices into $M_m$ of $m\times m$ matrices, the
{\sl Choi matrix} $\choi_\phi$ is given by
\begin{equation}\label{choi}
\choi_\phi=\sum_{i,j=1}^n e_{ij}\ot \phi(e_{ij})\in M_n\ot M_m,
\end{equation}
with the standard matrix units $\{e_{ij}\}$ of $M_n$.
This had been defined in 1967 by de Pillis \cite{dePillis} who had shown that $\phi$ preserves Hermitianity
if and only if $\choi_\phi$ is Hermitian.
During the seventies, two results appeared which are now fundamental tools in current
quantum information theory; Jamio\l kowski \cite{jam_72} gave a condition on $\choi_\phi$ under which $\phi$ is a positive linear map,
which is now called {\sl block positivity}; Choi \cite{choi75-10} showed that $\phi$ is completely positive if and only if $\choi_\phi$
is positive (semi-definite). The isomorphism $\phi\mapsto \choi_\phi$ is now called the Jamio\l kowski--Choi isomorphism.

Paulsen and Shultz \cite{Paulsen_Shultz} asked what happens when we replace matrix units in (\ref{choi})
by another basis $\calb=\{b_{ij}\}$ of $M_n$.
More precisely, they considered
\begin{equation}\label{p-s}
\choi^\calb_\phi=\sum_{i,j=1}^{n} b_{ij}\ot \phi(b_{ij})=(\id_n\ot\phi)\left(\textstyle\sum_{i,j} b_{ij}\ot b_{ij} \right)\in M_n\ot M_m
\end{equation}
and asked when the following statement
\begin{equation}\label{main}
\phi:M_n\to M_m\ {\text{\rm is completely positive}}\ \Longleftrightarrow\ \choi^\calb_\phi\in M_n\ot M_m\ {\text{\rm is positive}}
\end{equation}
is still true. For this purpose, they considered the map $C_\calb$ which sends matrix units to $\calb$;
the $n^2\times n^2$ matrix $[C_\calb]$ which represents the map $C_\calb$ with respect to the matrix units with a fixed order;
the transpose $[C_\calb]^T$ of the matrix $[C_\calb]$;
and finally the linear map $M_\calb:M_n\to M_n$ given by the $n^2\times n^2$ matrix $[C_\calb][C_\calb]^T$.
For a given $s\in M_n$, we define the map
$\ad_s:M_n\to M_n$ by
$$
\ad_s(x)=s^*xs,\qquad x\in M_n.
$$
Paulsen and Shultz \cite{Paulsen_Shultz} showed that if $M_\calb=\ad_s$ for some $s\in M_n$ then the statement (\ref{main}) holds
for every $m=1,2,\dots$.
We show that the converse is also true.

We note that the Choi matrix may be written by
\begin{equation}\label{k}
\choi_\phi=(\id_n\ot\phi)\left(\textstyle\sum_{i,j}e_{ij}\ot e_{ij}\right)=(\id_n\ot\phi)(|\omega\ran\lan\omega|),
\end{equation}
where $|\omega\ran=\sum_{i=1}^n |i\ran|i\ran\in\mathbb C^n\ot\mathbb C^n$.
Very recently, the author \cite{kye_comp_tensor}
considered the question what happens when we replace $|\omega\ran$ in (\ref{k}) by
another vector. He began with a nonsingular matrix $s=\sum_{i,j=1}^n s_{ij}{|i\ran\lan j|}\in M_n$
to associate the vector $\lan\tilde s|=\sum_{i,j=1}^n s_{ij} \lan i|\lan j|\in \mathbb C^n\ot \mathbb C^n$,
and defined
\begin{equation}\label{kkk}
\choi^s_\phi=(\id_n\ot\phi)(|\tilde s\ran\lan\tilde s|),
\end{equation}
to show that $\phi$ is completely positive if and only if $\choi^s_\phi$ is positive.
We show that non-singularity of $s$ is necessary to retain the correspondence between completely positivity of $\phi$
and positivity of $\choi^s_\phi$.

More generally, we will consider the question what happens when we replace
the matrix $\sum_{i,j}e_{ij}\ot e_{ij}$ in (\ref{k}) by
another matrix $\Sigma\in M_n\ot M_n$, to define
\begin{equation}\label{kk}
\choi^\Sigma_\phi=(\id_n\ot\phi)(\Sigma)\in M_n\ot M_n.
\end{equation}
We abuse notations in (\ref{p-s}), (\ref{kkk}) and (\ref{kk}).
For a given basis $\calb=\{b_{ij}\}$, we write $\Sigma_\calb=\sum_{i,j} b_{ij}\ot b_{ij}$ then
$\choi^\calb_\phi$ in (\ref{p-s}) is nothing but $\choi^{\Sigma_\calb}_\phi$ defined in (\ref{kk}).
The similar comment works for (\ref{kkk}) and (\ref{kk}).

If we take the linear map $\sigma:M_n\to M_n$ with
$\Sigma=\choi_\sigma$ then we have the relation
\begin{equation}\label{rel}
\choi^\Sigma_\phi=(\id_n\ot\phi)(\choi_\sigma)=\choi_{\phi\circ\sigma},
\end{equation}
which plays an important role in quantum information theory. See \cite[Section 4]{gks}.
Now, we will say that $\Sigma\in M_n\ot M_n$ {\sl satisfies the Choi correspondence} when
complete positivity of a map $\phi:M_n\to M_m$ is equivalent to positivity of $\choi^\Sigma_\phi\in M_n\ot M_m$
for every $m=1,2,\dots.$
We see that $\Sigma=\choi_\sigma$ satisfies the Choi correspondence if and only if
the following
\begin{equation}\label{c1}
\phi\in\cp[M_n,M_m]\ \Longleftrightarrow\ \phi\circ\sigma\in\cp[M_n,M_m]
\end{equation}
holds for every $m=1,2,\dots$ by (\ref{rel}), where $\cp=\cp[M_n,M_m]$ denotes the convex cone of all completely positive maps
from $M_n$ into $M_m$.

Following \cite{gks}, we use the bilinear pairing $\lan X,Y\ran=\tr(XY^\ttt)$ for matrices, and define
the bilinear pairing $\lan\phi,\psi\ran$ of two linear maps $\phi,\psi:M_n\to M_m$ by
$\lan\phi,\psi\ran=\lan\choi_\phi,\choi_\psi\ran$. The map $\phi^*:M_m\to M_n$ is
also defined by $\lan \phi^*(x),y\ran=\lan x,\phi(y)\ran$.
For a closed convex cone $K$ of linear maps, we also define $K^\circ$ as the convex cone consisting of all linear maps $\phi$
satisfying $\lan\phi,\psi\ran\ge 0$ for every $\psi\in K$. It is well known that $\cp^\circ=\cp$.
Then, we see that $\phi\circ\sigma\in\cp$ holds in (\ref{c1})
if and only if $\lan\phi\circ\sigma,\tau\ran\ge 0$ for every $\tau\in\cp$
if and only if $\lan\phi,\tau\circ\sigma^*\ran\ge 0$ for every $\tau\in\cp$ if and only if
$\phi\in(\cp\circ\{\sigma^*\})^\circ$.
Since $\cp^\circ=\cp$, we see that (\ref{c1}) holds if and only if
$\cp=\cp\circ\{\sigma^*\}$ if and only if the following
\begin{equation}\label{c2}
\cp[M_m,M_n]=\{\sigma\}\circ\cp[M_m,M_n]
\end{equation}
holds for every $m=1,2,\dots$.
Recall that a linear map $\sigma:M_n\to M_n$ is called a {\sl complete order
isomorphism} if $\sigma$ is a bijection and both $\sigma$ and
$\sigma^{-1}$ are completely positive.
When $m=n$, the identity map belongs to $\cp[M_n,M_n]$ in the condition (\ref{c2}),
and so we conclude that $\sigma$ is a complete order isomorphism.
Conversely, if $\sigma$ is a complete order isomorphism then
(\ref{c2}) holds for every $m=1,2,\dots$.
A typical example of complete order isomorphism is the map
$\ad_s:M_A\to M_A$ for a nonsingular $s\in M_A$, with the inverse map
$\ad_{s^{-1}}$.

\smallskip\noindent
{\bf Theorem 1.}
{\it
Suppose that $\Sigma\in M_n\ot M_n$ with $\Sigma=\choi_\sigma$ for a map $\sigma:M_n\to M_n$. Then the following are equivalent:
\begin{enumerate}
\item[(i)]
$\Sigma$ satisfies the Choi correspondence,
%that is, a linear map $\phi:M_n\to M_m$ is completely positive if and only if $\choi^\Sigma_\phi$ is positive for every $m=1,2,\dots$,
\item[(ii)]
$\phi\in\cp[M_n,M_m]$ if and only if $\phi\circ\sigma\in \cp[M_n,M_m]$ for every $m=1,2,\dots$,
\item[(iii)]
$\cp[M_m,M_n]=\{\sigma\}\circ\cp[M_m,M_n]$ for every $m=1,2,\dots$,
\item[(iv)]
$\sigma$ is a complete order isomorphism,
\item[(v)]
$\sigma=\ad_{s}$ for a nonsingular matrix $s\in M_n$,
\item[(vi)]
$\Sigma$ is a positive rank one matrix whose range vector has the full Schmidt rank,
\item[(vii)]
for any basis $\{|\xi_i\ran\}$ of $\mathbb C^n$, there exists a basis $\{|\eta_i\ran\}$ of $\mathbb C^n$ such that
$\Sigma=\sum_{i,j=1}^n |\xi_i\ran\lan\xi_j|\ot |\eta_i\ran\lan\eta_j|$,
\item[(viii)]
$\Sigma=\sum_{i,j=1}^n |\xi_i\ran\lan\xi_j|\ot |\eta_i\ran\lan\eta_j|$ for some bases $\{|\xi_i\ran\}$ and $\{|\eta_i\ran\}$ of $\mathbb C^n$.
\end{enumerate}
}
\smallskip

In \cite{kye_comp_tensor}, we have considered the case when
$\Sigma$ is a positive rank one matrix with the range vector $|\tilde s\ran\in\mathbb C^n\ot\mathbb C^n$, and showed that
if $s$ is non-singular then $\Sigma$ satisfies the Choi correspondence. The above theorem tells us that the converse also holds.

So far, we have seen that the statements (i), (ii), (iii) and (iv) are equivalent. It was shown in
\cite{Paulsen_Shultz} that (iv) implies (v). We show here in another way.
Suppose that
$\sigma=\sum_{i\in I}\ad_{s_i}$ is a complete order isomorphism with the
completely positive inverse $\sigma^{-1}=\sum_{j\in J} \ad_{t_j}$. Then we
have
$\sum_{i,j}\ad_{t_js_i}= \sum_{i,j} \ad_{s_i}\circ\ad_{t_j} =\id$.
Because the identity map generates an extreme ray of the convex cone $\cp[M_n,M_n]$,
there exist positive numbers $\alpha_{(i,j)}$ such that
$\alpha_{(i,j)}\ad_{t_js_i}=\id$ for $(i,j)\in I\times J$,
and so, we have
$\sqrt{\alpha_{(i,j)}}z_{(i,j)}{t_js_i}=I$ for every $(i,j)\in I\times J$, with complex numbers $z_{(i,j)}$ of modulus one.
Therefore, we see that all $s_i$ coincide up to scalar multiplications.
Furthermore, the invertibility of $\ad_s$ implies non-singularity of $s$. This shows that (iv) and (v) are equivalent.
By the relation $\choi_{\ad_s}=|\tilde s\ran\lan\tilde s|$ in \cite{gks}, we also have (v) $\Longleftrightarrow$ (vi).
For a given basis $\{|\xi_i\ran\}$ of $\mathbb C^n$, every vector in $\mathbb C^n\ot\mathbb C^n$
is of the form $\sum_{i=1}^n|\xi_i\ran\ot |\eta_i\ran$ for a family $\{|\eta_i\ran\}$ of vectors, which
has the full Schmidt rank if and only if  $\{|\eta_i\ran\}$ is
a basis. This shows the direction (vi) $\Longrightarrow$ (vii). The remaining implications (vii) $\Longrightarrow$ (viii) $\Longrightarrow$ (vi)
are clear.

We note that $\phi\mapsto \choi^\Sigma_\phi$ is a bijection if and only if $\phi\mapsto\phi\circ\sigma$
is a bijection between the linear mapping space $L(M_n,M_m)$
if and only if the map $\sigma\in L(M_n,M_n)$ is invertible. When $\sigma=\ad_s$, we also see that $\sigma$ is invertible if and only if
$s$ is nonsingular. Therefore, we see that the Choi correspondence of $\Sigma$ implies that
$\phi\mapsto \choi^\Sigma_\phi$ is a bijection.

Considering the condition in \cite{Paulsen_Shultz},
we suppose that $\calb=\{b_j:j\in J\}$ is an ordered family of $n\times n$ matrices with the cardinality $n^2$, and
define the linear map $C_\calb: M_n\to M_n$ by
$b_j=C_\calb (e_j)=\sum_{i\in J}x_{ij}e_i$ for $j\in J$,
with an ordered matrix unit $\{e_j:j\in J\}$.
By the relation $C_\calb^T(e_j)=\sum_{i\in J}x_{ji}e_i$, we have
$$
M_\calb(e_j)=C_\calb C_\calb^T(e_j)=\sum_{i\in J} x_{ji} C_\calb(e_i)
=\sum_{i\in J}x_{ji}\sum_{k\in J} x_{ki}e_k.
$$
Therefore, it follows that
$$
\begin{aligned}
\sum_{j\in J} e_j\ot M_\calb(e_j)
&=\sum_{j\in J}\sum_{i\in J}\sum_{k\in J} x_{ji}x_{ki} e_j\ot e_k\\
&=\sum_{i\in J}\left(\sum_{j\in J} x_{ji}e_j\right)\ot \left(\sum_{k\in J} x_{ki}e_k\right)
=\sum_{i\in J} b_i\ot b_i.
\end{aligned}
$$

Let $\calb=\{b_{ij}:i,j=1,2,\dots,n\}$ be a family of matrices in $M_n$. If we put
$$
\Sigma=\sum_{i,j=1}^n b_{ij}\ot b_{ij}=\sum_{i,j=1}^n e_{ij}\ot M_\calb(e_{ij})=\choi_{M_\calb},
$$
then $\choi^\Sigma_\phi$ in (\ref{kk}) is nothing but $\choi^\calb_\phi$ defined in (\ref{p-s}).
Therefore, the statement (\ref{main}) holds for every $m=1,2,\dots$ if and only if $M_\calb=\ad_s$ for a nonsingular $s\in M_n$
if and only if $\sum_{i,j=1}^n b_{ij}\ot b_{ij}=\sum_{i,j=1}^n |\xi_i\ran\lan\xi_j|\ot |\eta_i\ran\lan\eta_j|$ for some bases
$\{|\xi_i\ran\}$ and $\{|\eta_i\ran\}$ if and only if there exists a basis $\{|\zeta_i\ran\}$ such that $b_{ij}=|\zeta_i\ran\lan\zeta_j|$
for each $i,j=1,2,\dots,n$. In case when $\{b_{ij}\}$ is a basis of $M_n$, the matrix $s$ in this discussion
is automatically non-singular, by the invertibility of $M_{\mathcal B}$. We summarize as follows:

\smallskip\noindent
{\bf Theorem 2.}
{\it
For a basis $\calb=\{b_{ij}\}$ of $M_n$, the following are equivalent:
\begin{enumerate}
\item[(i)]
$\phi: M_n\to M_m$ is completely positive if and only if $\choi^\calb_\phi\in M_n\ot M_m$ is positive,
\item[(ii)]
$M_\calb=\ad_s$ for a nonsingular matrix $s\in M_n$,
\item[(iii)]
$M_\calb=\ad_s$ for a matrix $s\in M_n$,
\item[(iv)]
there exists a basis $\{|\zeta_i\ran\}$ of $\mathbb C^n$ such that $b_{ij}=|\zeta_i\ran\lan\zeta_j|$ for $i,j=1,2,\dots,n$.
\end{enumerate}
}
The author is grateful to Kyung Hoon Han for helpful coments on the draft.

\end{document}